\documentclass[10pt,twoside,leqno]{article}

\usepackage[english]{babel}
\usepackage[utf8]{inputenc}
\usepackage[pdfpagelabels=true]{hyperref}

\def\twoupletter{%
\usepackage[landscape]{geometry}%
\usepackage{2up}
\setlength{\paperwidth}{11.0in}
\setlength{\paperheight}{8.5in}
\setlength{\topmargin}{-0.75in}
\setlength{\headsep}{0.15in}
\setlength{\footskip}{0.3in}
\setlength{\textheight}{7.2in}
\setlength{\oddsidemargin}{-0.65in}   
\setlength{\evensidemargin}{-0.8in}
\setlength{\textwidth}{5.0in}   
\source{\magstep0}{5.5in}{8.267in}
\target{\magstep0}{11.0in}{8.5in}
\newcounter{dpage}
\def\oddno##1{\setcounter{dpage}{##1}\multiply \value{dpage}  by 2 \addtocounter{dpage}{-1}\thedpage}
\def\evenno##1{\setcounter{dpage}{##1}\multiply \value{dpage}  by 2 \addtocounter{dpage}{0}\thedpage}}

\def\twoupa4{%
\usepackage[landscape]{geometry}
\usepackage{2up}
\setlength{\paperwidth}{11.692in}
\setlength{\paperheight}{8.267in}
\setlength{\topmargin}{-0.75in}
\setlength{\headsep}{0.15in}
\setlength{\textheight}{7.2in}
\setlength{\oddsidemargin}{-0.78in} 
\setlength{\evensidemargin}{-0.78in} 
\setlength{\textwidth}{5.1in}  
\source{\magstep0}{5.5in}{8.267in}
\target{\magstep0}{11.692in}{8.267in}
\newcounter{dpage}
\def\oddno##1{\setcounter{dpage}{##1}\multiply \value{dpage}  by 2 \addtocounter{dpage}{-1}\thedpage}
\def\evenno##1{\setcounter{dpage}{##1}\multiply \value{dpage}  by 2 \addtocounter{dpage}{0}\thedpage}}

\twoupletter

\usepackage{amsmath,amssymb,fancyhdr}

\DeclareMathOperator{\sing}{sing}
\DeclareMathOperator{\reg}{reg}
\DeclareMathOperator{\dvg}{div}

\DeclareMathOperator{\graph}{graph}
\DeclareMathOperator{\spt}{spt}
\DeclareMathOperator{\density}{density}

\def\widebar#1{\vbox{\hrule height.1pt\ialign{##\crcr\hrulefill\crcr\noalign{\kern1.6pt\nointerlineskip}
 $\hfil\displaystyle{#1}\hfil$\crcr}}} %
\def\swidebar#1{\vbox{\scriptsize\hrule height.1pt\ialign{##\crcr\hrulefill\crcr\noalign{\kern1.0pt\nointerlineskip}
 $\hfil\displaystyle{#1}\hfil$\crcr}}} %

\newcounter{sequation}[section] 
\renewcommand{\thesequation} {\arabic{section}.\arabic{sequation}} 
\def\dl#1{\refstepcounter{sequation}\label{#1}\leqno\textrm{\thesequation}} 
\newcounter{pequation} 
\def\dtg#1{\refstepcounter{sequation}\label{#1}\tag*{\textrm{\thesequation}}} 

\newenvironment{state}[1]{\par\smallbreak\noindent {\bf#1} %
\textit\bgroup}{\egroup \par\ifdim\lastskip<\medskipamount \removelastskip\penalty55\medskip\fi} %
 
\def\tfrac#1#2{\textstyle\frac{#1}{#2}}

\def\S{$\mathsection$}

\def\R{\mathbb{R}}\def\C{\mathbb{C}}\def\Sph{\mathbb{S}}
\renewcommand{\epsilon}{\varepsilon}

\begin{document}

\title{\vskip-.6in           Singularities of minimal submanifolds }  \author{\scshape                         Leon Simon }
                                                \date{\vspace{-.0in}}%

\maketitle

\thispagestyle{empty}

\begin{abstract}
\noindent
After quick survey of some key results and open questions about the structure of singularities of minimal surfaces, we discuss recent
work~\cite{Sim23} on singularities of stable minimal hypersurfaces, including some simplifications of the main technical discussion
in~\cite{Sim23}.
\end{abstract}

\section{Introductory Remarks and Main Theorem}

The most general notion (at least in the natural setting where multiplicities are assumed to be integer-valued) of an
$n$-dimensional minimal submanifold of $U$, $U$ open in $\R^{n+k}$, is the ``integer multiplicity stationary
varifold'' i.e.\ a countably $n$-rectifiable set $M$ of locally finite $n$-dimensional Hausdorff measure in $U$, equipped
with a locally bounded density function $\theta$ which is positive integer valued $\mathcal{H}^{n}$-a.e.\ on $M$ (and
$0$ on $U\setminus M$) such that the first variation of the area $\mu(M)=\int_{M}d\mu$ vanishes in $U$, where
$d\mu=\theta d\mathcal{H}^{n}$; thus
\[%
\frac{d}{dt}_{\bigl|t=0}\int_{\Psi_{t}(M)}\theta\circ\Psi_{t}^{-1}(x)\,d\mathcal{H}^{n} =0, %
\dl{1st-var}
\]%
where $\Psi_{t}(x)=x+tX(x)$, where $X=(X^{1},\ldots,X^{n+k})$ is any $C^{1}$ vector field on $U$ with compact
support in $U$. Via a straightforward computation using the area formula,~\ref{1st-var} is equivalent to the first variation
identity
\[%
\int_{M} \dvg_{M}X\, d\mu=0, %
\dl{1st-formula}
\]%
where $\dvg_{M}X$ is the tangential divergence, computed at points $x\in M$ by the formula
$\dim_{M}X=\sum_{i=1}^{n}\tau_{i}\cdot D_{\tau_{i}}X$, where $\tau_{1},\ldots,\tau_{n}$ is any orthonormal basis of
$T_{x}M$, the approximate tangent space of $M$ at $x$ (which exist for $\mathcal{H}^{n}$-a.e.\ $x\in M$).

In view of the well-known monotonicity of $(\omega_{n}\rho^{n})^{-1}\mu(B_{\rho}(x))$ ($\omega_{n}=$ volume of
unit ball in $\R^{n}$), for each $x\in U$ we can replace $\theta$ with the preferred representative
$\Theta^{n}(x)=\lim_{\rho\downarrow 0}(\omega_{n}\rho^{n})^{-1}\mu(B_{\rho}(x))\,(=\theta(x)$ for
$\mathcal{H}^{n}$-a.e.\ $x$) and replace $M$ by $\{x\in U:\Theta^{n}(x)\ge 1\}$, which is a closed set because
$\Theta^{n}$ is upper semi-continuous on $U$, again by the monotonicity of $\rho^{-n}\mu(B_{\rho}(x))$.

So the general integer multiplicity stationary setting in an open set $U$ reduces to the case when $M$ is a closed
rectifiable subset of $U$ with locally finite $\mathcal{H}^{n}$-measure such that $M$, with multiplicity
$\Theta^{n}(x)=\lim_{\rho\downarrow 0}(\omega_{n}\rho^{n})^{-1}\smash{\int_{M\cap
B_{\rho}(x)}\theta\,d\mathcal{H}^{n}}$ (which exists everywhere in $U$ and $=0$ on $U\setminus M$), $\Theta^{n}$ is
upper semi-continuous, and the area functional (with measure $d\mu=\Theta^{n}d\mathcal{H}^{n}$ on $M$) is stationary
with respect to compactly supported deformations of the identity in the ambient space $U$.

As usual we define
\[%
\quad\reg M=\text{regular set}, %
\dl{def-reg} %
\]%
i.e.\ the set of $x\in M$ such that $M\cap B_{\rho}(x)$ is a smooth embedded submanifold for some $\rho>0$, and %
\[%
\sing M=M\setminus \reg M,  %
\dl{def-sing}
\]%
so $\sing M$ is by definition a closed subset of $U$.

In this general setting not a lot is known about the singular set.  Intuitively one would perhaps expect that $\sing M$
should have $n$-dimensional measure zero, but remarkably such a basic question has remained open since Almgren's initial
development in the 1960's of the theory of varifolds, a theory later fully developed by Allard in his famous 1972
paper~\cite{All72}.

The best that is known in this general setting is that
\[%
\sing M \text{ is nowhere dense in $M$, } %
\dl{ndense}
\]%
i.e.\ $\reg M$ is dense in $M$, which is clear, because the Allard regularity theorem guarantees that any point where
$\Theta^{n}$ is a local minimum must be a regular point---and every open ball intersecting $M$ contains points where
$\Theta^{n}$ ($=$ integer a.e.) is a minimum in that ball. Thus
\[%
M=\widebar{\reg M} %
\dl{m0}
\]%
and (by the constancy theorem) the density function $\Theta^{n}$ is a constant positive integer on each connected
component of $\reg M$.

The situation is much better for special classes of minimal submanifolds $M$:

Due to pioneering work of De Giorgi, Reifenberg, Federer, Almgren, and James Simons, it is known that codimension 1
minimizing $M\subset\R^{n+1}$ cannot have singularities if dimension $M\le 6$, singularities are discrete for $n=7$, and
the dimension of the singular set is $\le n-7$ for $n\ge 8$, and in codimension $\ge 2$ the singular set of mod-2
minimizers is discrete in case $n=2$ and has dimension $\le n-2$ in case $n\ge 3$.

The above results for codimension 1 minimizers were extended to stable $M$ in~\cite{SchSY75} and~\cite{SchS81}.

It was shown in \cite{Sim93}, \cite{Sim95b} that in all the above cases the singular set is rectifiable.

The case of oriented minimizing $M$ (as distinct from mod 2 minimizers) in codimension $\ge 2$ is much more difficult
(due to branching), and Almgren in the early 1980's finally proved that the singular set of $M$ has dimension $\le n-2$.
Almgren's proof was simplified and shortened by De Lellis-Spadaro~\cite{DelS14a}, \cite{DelS14b}, \cite{DelS14c}. 
Finally~Krummel-Wickramasekera and (independently) De Lellis-Minter-Skorobogatova in 2023 proved that the singular
set is $(n-2)$-rectifiable.

Thus for all the special classes mentioned above, up to a set of measure zero in the appropriate dimension, the singular set
$K=\sing M$ can be decomposed into a countable union of pieces, each of which is contained in an embedded $C^{1}$
submanifold of the ambient space.

The obvious question then is: what do these pieces look like? e.g.\ do they typically have positive measure in the $C^{1}$
submanifold, or can they be fractional dimensional?

In high enough dimensions it has recently been established there are in fact examples where the singular set can be very bad
(e.g.\ fractional dimensional): indeed with respect to a suitable smooth metric for $\R^{n+1+\ell}$, $n\ge 7,\,\ell\ge 1$,
there are complete stable minimal hypersurfaces $M\subset\R^{n+1+\ell}$ with $\sing M=K$, $K$ any preassigned closed
subset of $\{0\}\times \R^{\ell}$; more precisely

\begin{state}{\bf{}Main Theorem}{\rm(\cite{Sim23}).}  Suppose $n\ge 7$, $\ell\ge 1$, $K\subset\R^{\ell}$ is an
arbitrary closed subset, and $\epsilon>0$. Then there is a smooth metric
$g=\smash{\sum_{i,j=1}^{n+1+\ell}g_{ij}dx^{i}dx^{j}}$ on $\R^{n+1+\ell}$ with $|g_{ij}-\delta_{ij}|<\epsilon$ and
$|D^{k}g_{ij}|<C_{k}\epsilon$ for each $k=1,2,\ldots$ and such that there is a smooth embedded hypersurface $M
\subset\R^{n+1+\ell}$ which is minimal and stable with respect to the ambient metric $g$ and which has $\sing
M=\{0\}\times K$.
\end{state}

Whether or not such phenomena can happen with respect to real analytic metrics $g$ (or the Euclidean metric) is open,
although G\'abor Sz\'ekelyhidi~\cite{Szk22} shows that all the non-isolated singularity examples of the type discussed here
cannot occur in case the metric is the standard Euclidean metric.

In lower dimensions (dimension $M\ge 3$ and codimension at least 3) Zhenhua Liu~\cite{Liu23} has proved using
calibration methods that there are examples of $M$ which are area minimizing with respect to a smooth metric and which
have fractional dimensional singular sets $K$ which have pieces representing the various ``Almgren strata'' in all the
appropriate dimensions, and $M$ minimizes relative to a suitable smooth metric. For any point $x\in K$, for small
enough $\rho>0$, $M\cap B_{\rho}(x)$ has the form $\Sigma_{1}\cup\Sigma_{2}$, $\text{\normalsize
sing\,}\Sigma_{1}=\text{\normalsize sing\,}\Sigma_{2}=\emptyset$, $\Sigma_{1}$, $\Sigma_{2}$ touch tangentially along
$K\cap B_{\rho}(x)=\Sigma_{1}\cap \Sigma_{2}$, so $\text{\normalsize sing\,}M\cap
B_{\rho}(x)=\Sigma_{1}\cap\Sigma_{2}\,(=K\cap B_{\rho}(x))$.

\smallskip

We are here going to outline the main details of the proof of the Main Theorem which was originally given
in~\cite{Sim23}.  For simplicity of notation in the remainder of the present discussion we take $\ell=1$, and in \S2 we
also make some simplifications of the existence discussion for solutions of the Symmetric Minimal Surface Equation
(SME), which is a key part of the proof in~\cite{Sim23}.

The proof to be outlined here uses the cylindrical cone
\[%
\C_{0}\times \R, \dl{def-cyl}
\]%
where $\C_{0}$ is the ``James Simons cone''
\[%
\C_{0}=\left\{(\xi,\eta)\in \R^{m+1}\times\R^{m+1} :|\eta|=|\xi|\right\}, %
\]%
In~\cite{Sim23} we allowed $\C_{0}$ to be any one of the cones
$\{(\xi,\eta)\in\R^{m+1}\times\R^{n+1}:n|\xi|^{2}=m|\eta|^{2}\}$ with $m,n\ge 2$ and $m+n\ge 6$, but the argument
outlined here for the special case of Simons cone involves only notational changes from this more general setting.

\section {\!\!\!Special Smoothing of $\sing \C_{0}\times\R$ in 
$\R^{2m+2}\times(\R\,\setminus K)$}

Initially we assume that $R>1$ (large), and the indicator function $\chi_{K}$ of $K$ is $2R$-periodic.  Eventually, at the
end of the discussion, we let $R\to \infty$ to handle the general case of arbitrary closed $K$.

Let $h:\R\to [0,\frac{1}{2}]$ be a $C^{\infty}$ function with $h>0$ on $\R\setminus K$ and $h=0$ on $K$, and $h$
periodic with period $2R$.

We can smooth out the singularities of $\C_{0}\times\R$ along $\{0\}\times(\R\setminus K)$ by using surfaces with cusps
at end points of the intervals of $\R\setminus K$; more specifically let $\Omega$ be the domain in $\R^{m+1}\times\R$
given by
\[%
\Omega=\{(\xi,y)\in\R^{m+1}\times\R:|\xi|<h(y)\}, %
\]%
and let $u$ be a non-negative $C^{0}(\R^{m+1}\times\R)$ function, $2R$-periodic in the $y$ variable, which is positive
and $C^{\infty}$ on $\bigl((\R^{m+1}\times(\R\setminus(\{0\}\times K)))$, and $u(\xi,y)=|\xi|$ on $(\R^{m+1}\times
\R)\setminus\Omega$. We then make an $m$-fold rotation of $\graph u$ to give a hypersurface in
$\R^{m+1}\times\R^{m+1}\times\R$; that is, we replace $\graph u$ by the \emph{symmetric graph} $SG(u)$, where
\[%
SG(u) = \{(\xi,\eta,y)\in\R^{m+1}\times\R^{m+1}\times\R:|\eta|=u(\xi,y)\}. %
\]%
$SG(u)$ is then a hypersurface in $\R^{m+1}\times\R^{m+1}\times\R$ with $\sing SG(u)=\{0\}\times\{0\}\times K$. 
We aim to choose $u$ so that $SG(u)$ is stationary (i.e.\ has zero first variation of area) with respect to some smooth
metric for $\R^{m+1}\times\R^{m+1}\times\R$ to be determined; for this procedure to succeed we need some special
properties (\ref{props-u}, \ref{ext-u} below) of $u$.

As a preliminary to this discussion we need to make some remarks about the possibility that a symmetric graph $SG(v)$,
where $v$ is positive and $C^{\infty}$ on an open set $U\subset\R^{m+1}\times\R$, can be a minimal hypersurface with
respect to the standard Euclidean metric for $\R^{m+1}\times\R^{m+1}\times\R$. Such minimality would require that
the area of $SG(v)$ (i.e.\ $\mathcal{H}^{2m+2}(SG(v))$) is stationary with respect to ambient compactly supported
deformations of $SG(v)$. By the area formula, applied to the map $\{(\xi,\omega,y):(\xi,y)\in U,\, \omega\in \Sph^{m}\}
\to \R^{m+1}\times\R^{m+1}\times\R$ which takes $(\xi,\omega,y)$ to $(\xi,v(\xi,y)\omega,y)$, we have
\[%
\mathcal{H}^{2m+2}(SG(v)) = \sigma_{m}\int_{U}\sqrt{1+|Dv|^{2}}\,v^{m}\,d\xi dy, %
\dl{area-sgv}
\]%
where $\sigma_{m}$ is the volume of the unit sphere $\Sph^{m}\subset\R^{m+1}$.  Minimality of $SG(v)$ is thus
equivalent to the requirement that $v$ satisfies the Euler-Lagrange equation of the functional on the right of~\ref{area-sgv}.

The Euler-Lagrange equation, which we call the Symmetric Minimal Surface Equation (SME), is in fact
$\mathcal{M}(v)=0$, where
\[%
\mathcal{M}(v)= \sqrt{1+|Dv|^{2}}\,%
\Bigl({\sum}_{i=1}^{m+1}D_{\xi^{i}}\bigl(\dfrac{D_{\xi^{i}}v}{\sqrt{1+|Dv|^{2}}}\bigr) %
+D_{y}\bigl(\dfrac{D_{y}v}{\sqrt{1+|Dv|^{2}}}\bigr)\Bigr) -\frac{m}{v}. %
\dl{SME}
\]%
Notice in particular this ensures
\[%
\Bigl({\sum}_{i=1}^{m+1}D_{\xi^{i}}\bigl(\dfrac{D_{\xi^{i}}v}{\sqrt{1+|Dv|^{2}}}\bigr) %
+D_{y}\bigl(\dfrac{D_{y}v}{\sqrt{1+|Dv|^{2}}}\bigr)\Bigr) > 0, %
\]%
for positive solutions of the SME, so such solutions satisfy the strict maximum principle; in particular $v$ cannot attain
an a local maximum at an interior point of the domain of $v$.

Notice that for functions $v(\xi,y)=\varphi(|\xi|)>0$ (i.e.\ positive functions which are independent of the $y$ variable and
which can be expressed as a function of the single variable $r=|\xi|$) the functional in~\ref{area-sgv}, assuming
$U=\{(\xi,y)\:|\xi|<\rho,\,|y|<\rho\}$, can be written
\[%
2\rho\sigma_{m}\int_{0}^{\rho}\sqrt{1+(\varphi')^{2}}\,\varphi^{\,m}r^{m}\,drdy, %
\dl{area-vph}
\]%
and in this case $\mathcal{M}(\varphi)=\frac{\varphi''(r)}{1+(\varphi'(r))^{2}} + \frac{m}{r}\varphi'(r) -
\frac{m}{\varphi(r)}$, so the SME is just the ODE
\[%
\frac{\varphi''(r)}{1+(\varphi'(r))^{2}} + \frac{m}{r}\varphi'(r) - \frac{m}{\varphi(r)}=0, %
\dl{vph-eqn}
\]%
and $SG(\varphi)$ is a cylinder (i.e.\ independent of the $y$-variable) lying on one side of the Simons cylinder
$\C_{0}\times\R$. As shown in~\cite{Sim23}, such a $\varphi$ is convex, has a positive minimum (which by scaling we
can, and we shall, assume to be $1$) at $r=0$, and is asymptotic to $r$ at $\infty$. In fact by some ODE manipulations
(as in~\cite{Sim23}),
\[%
\,\left\{\begin{aligned}%
&\,\, \varphi''(r) >0\,\,\forall r\ge 0,\,\,\, r<\varphi(r)<1+ r,\, 0<\varphi'(r)<1\,\,\,\forall r>0 \\ %
&\,\,\varphi(r)- r\sim \kappa r^{-\gamma}, \,\varphi(r)-r\varphi'(r)\sim %
\kappa(1+\gamma) r^{-\gamma}\!,\,\varphi''(r)\sim \kappa\gamma(\gamma+1)r^{-\gamma-2} %
\end{aligned}\right.%
\dl{vph-props}
\]%
as $r\to \infty$, where $\kappa=\kappa(m)$ is a positive constant and
\[%
\gamma=(m-\tfrac{1}{2}) -\sqrt{(m-\tfrac{1}{2})^{2}-2m}\,\, \text{ ($=2$ when $m=3$, $>1$ when $m>3$).} %
\]%
Also, since $\frac{d}{dr}(\varphi(r)-r\varphi\hskip1pt'(r))=-r\varphi\hskip1pt''(r)<0$, we see that $\varphi(r)
-r\varphi\hskip1pt'(r)$ is strictly decreasing and hence, since $\lim_{r\to\infty}(\varphi-r\varphi\hskip1pt')=0$
by~\ref{vph-props},
\[%
0<\varphi(r)-r\varphi\hskip1pt'(r)<\varphi(0)=1 \text{ for all } r>0.  %
\dl{b}
\]%
Also, all the geometric rescalings
\[%
\varphi_{\tau}(r)=\tau\varphi(r/\tau),\quad \tau >0, %
\dl{vph-tau}
\]%
also satisfy~\ref{vph-eqn}, and $\varphi_{\tau}(r)>r$, $\varphi_{\tau}(r)\downarrow \varphi_{0}(r)$ as $\tau\downarrow
0$, where $\varphi_{0}(r)=r$, which has symmetric graph equal to the Simons cylinder $\C_{0}\times \R$. Thus the
family $SG(\varphi_{\tau}),\,\tau>0$, foliates all of the volume $|\eta|>|\xi|$ on one side of the Simons cylinder. Also
\[%
\varphi_{\tau}(r) - r<\tau, %
\dl{vtau-r}
\]%
because $\varphi_{\tau}(r) - r=\tau\varphi(r/\tau)-r <\tau(r/\tau+1)-r=\tau$ by the second inequality in~\ref{vph-props}.

Using the fact that the Simons cone $\C_{0}$ is strictly stable and strictly minimizing (see \cite{HarS85} for discussion),
one can check, as in~\cite{Sim23}, that $SG(\varphi)$ and the rescalings $SG(\varphi_{\tau})$ are strictly stable and hence
also the cylinders
\[%
SG(\varphi)\times\R \text{ and the rescalings } SG(\varphi_{\tau})\times\R \text{ are strictly stable, } %
\dl{stab-vph}
\]%
meaning that there is a fixed $\lambda>0$ such that the first eigenvalue of the linearization at $u=\varphi_{\tau}$ of the
SME operator $\mathcal{M}(u)$ on any bounded domain in $\R^{m+1}\times\R$ (with respect to zero Dirichlet
boundary data) is $\ge \lambda$ (independent of $\tau$).

By applying all of the above discussion to solutions $\widetilde{\varphi}$ of the Euler-Lagrange equation of the functional
\[%
\!\int_{0}^{\rho}\sqrt{1+(\widetilde{\varphi}\hskip1pt')^{2}}\, %
\widetilde{\varphi}^{\,\tilde m}r^{\tilde{m}}\,drdy %
\dl{til-fnl}
\]%
(instead of the original area functional in~\ref{area-vph}), where $\widetilde m=m/(1+\eta)$ with $\eta>0$ small, we
obtain in place of $\varphi$ the function $\widetilde{\varphi}$ which satisfies all the properties in~\ref{vph-props} with
$\widetilde{m}$ in place of $m$ and with $\widetilde\gamma$ in place of $\gamma$, where
\[%
\widetilde\gamma=(\widetilde m-\tfrac{1}{2}) -\sqrt{(\widetilde m-\tfrac{1}{2})^{2}-2\widetilde m}\,\,(>\gamma),%
\]%
and $\widetilde\varphi$ also satisfies the Euler-Lagrange equation of the functional in~\ref{til-fnl} and hence
\[%
(1+\eta)\frac{\widetilde{\varphi}\hskip1pt''(r)}{1+(\widetilde{\varphi}\hskip1pt'(r))^{2}} %
+ \frac{m}{r}\widetilde{\varphi}\hskip1pt'(r) - \frac{m}{\widetilde{\varphi}(r)}=0, %
\]%
giving
\[%
\mathcal{M}(\widetilde\varphi) =-\frac{\eta}{1+(\widetilde{\varphi}\hskip1pt')^{2}} %
\,\widetilde\varphi\,''< -\frac{\eta}{2}\,\widetilde\varphi\,'' <0.  %
\dl{vph-super}
\]%
Notice also that with $\epsilon\in(0, \epsilon_{0})$, where $\epsilon_{0}=\epsilon_{0}(m)\in (0,1)$ is sufficiently small,
we have
\[%
\varphi_{\epsilon^{\alpha}}(r)\le \widetilde{\varphi}_{\epsilon}(r),\quad \forall r\le 1, %
\dl{pwr-a}
\]%
provided we take $\alpha>\frac{\tilde\gamma+1}{\gamma+1}\,(>1)$. Since we can take $\widetilde{m}$ as close to $m$ as
we please, this effectively means we can take any $\alpha>1$.

The inequality~\ref{vph-super} enables us the construct a family of supersolutions of the SME on any of the domains
$\Omega_{\epsilon}$,
\[%
\Omega_{\epsilon}=\{(\xi,y):|\xi|<h(y)+\sqrt{\epsilon}\},\quad\epsilon\in [0,\tfrac{1}{2}\big].  %
\dl{om-d}
\]%
For any $v=v(r,y)$ on $\Omega_{\epsilon}$ ($r=|\xi|$),
\begin{align*}%
\mathcal{M}(v) %
\dtg{Mv}&=\frac{1+v_{y}^{2}} {1+v_{r}^{2}+v_{y}^{2}}v_{rr}+ \frac{1+v_{r}^{2}} {1+v_{r}^{2}+v_{y}^{2}}v_{yy} %
-2\frac{v_{r}v_{y}v_{ry}} {1+v_{r}^{2}+v_{y}^{2}}+m\Big(\frac{v_{r}} {r}-\frac{1}{v}\Big) \\ %
&=\frac{1} {1+v_{r}^{2}}v_{rr} +m\Big(\frac{v_{r}} {r}-\frac{1}{v}\Big) + %
\frac{v_{r}^{2}v_{y}^{2}} {(1+v_{r}^{2}+v_{y}^{2})(1+v_{r}^{2})} v_{rr} \\ %
\noalign{\vskip-4pt}
&\hskip1.8in + \frac{1+v_{r}^{2}} {1+v_{r}^{2}+v_{y}^{2}}v_{yy}-2\frac{v_{r}v_{y}v_{ry}} {1+v_{r}^{2}+v_{y}^{2}} \\ %
&\le \frac{1} {1+v_{r}^{2}}v_{rr} +m\Big(\frac{v_{r}} {r}-\frac{1}{v}\Big) + |v_{y}||v_{rr}|+ |v_{yy}| +|v_{ry}|.
\end{align*}%
We select $v$ to have the special form
\[%
v(r,y)=\psi(y))\widetilde\varphi\bigl(r/\psi(y)\big) \dl{v-choice}
\]%
with $\psi>0$ on $\R$; then
\[%
\begin{aligned}%
&v_{r}=\widetilde{\varphi}\hskip1pt'(r/\psi),\,\, v_{y} =\psi' b, \,\,\, %
v_{yy}=(r^{2}(\psi')^{2}/\psi^{2})\psi^{-1}\widetilde{\varphi}\hskip1pt''(r/\psi)+\psi'' b, \\   %
&v_{rr}=\psi^{-1}\widetilde{\varphi}\hskip1pt''(r/\psi),\,\,  %
v_{ry}=-(r\psi'/\psi)\psi^{-1}\widetilde{\varphi}\hskip1pt''(r/\psi),\,\,  %
\end{aligned}%
\]%
where $b= \widetilde{\varphi}(r/\psi)-(r/\psi)\widetilde{\varphi}\hskip1pt'(r/\psi)
=(\widetilde{\varphi}(t)-t\widetilde{\varphi}\hskip1pt'(t))|_{t=r/\psi}$.

Evidently $0<b<1$ by~\ref{b}, and since $\widetilde{\varphi}(t)-t\widetilde{\varphi}\hskip1pt'(t)\le
b_{0}(1+t^{2})\widetilde{\varphi}\hskip1pt''(t)$ for suitable constant $b_{0}$ by~\ref{vph-props},  we also have $b\le
b_{0}(1+r^{2}/\psi^{2})\widetilde{\varphi}\hskip1pt''(r/\psi)$. Using these facts together with~\ref{vph-super} in~\ref{Mv},
with the choice of $v$ in~\ref{v-choice}, we obtain
\[%
\mathcal{M}(v) \le\left(\!-\tfrac{\eta}{2} + |\psi'|+r^{2}(\psi')^{2}/\psi^{2} %
+ b_{0}(1+r^{2}/\psi^{2})\psi|\psi''| +(r|\psi'|/\psi) \right) \psi^{-1}\widetilde{\varphi}\hskip1pt''(r/\psi). %
\dl{Mv-bd}
\]%
Henceforth we make the additional assumptions that, with $\beta\in (0,2^{-7}]$ to be chosen (small), $\sqrt{h}\in
C^{\infty}$ and
\[%
\sqrt{\epsilon}+|\sqrt{h}|_{C^{1}}+|h|_{C^{2}} \le \beta, 
\quad |D^{k}h| \le C_{k}\beta, \,\, k\ge 3
\dl{h-beta}
\]%
(which of course we can arrange by simply taking a new function $h$ equal to the square of the original function and then
multiply by a suitable factor), and select $\psi=\psi_{s}$, where $s\ge \epsilon$ and
\[%
\psi_{s}(y)
=\begin{cases}%
\,\,\, s+ \beta e^{-(\sqrt{s}+h(y))^{-1/2}}, & s<1 \\  %
\sqrt{s}+ \beta e^{-(\sqrt{s}+h(y))^{-1/2}}, & s\ge 1. %
\end{cases}%
\dl{psi-d}
\]%
Also let $v_{s}$ be the corresponding $v$; that is,
\[%
v_{s}(r,y)=\psi_{s}(y)\widetilde{\varphi}(r/\psi_{s}(y)).   %
\dl{def-ved}
\]%
Notice that by~\ref{b} (with $\widetilde{\varphi}$ in place of $\varphi$) this is an increasing function of $s$ so in
particular 
\[%
v_{s}\ge v_{\epsilon}, \quad s\ge \epsilon.
\dl{mon-v}
\]%
With $\psi=\psi_{s}$ as in~\ref{psi-d}, we have
\[%
\psi'' = \beta\Big(\big(-\tfrac{3}{4}(h+\sqrt{s})^{-5/2}+\tfrac{1}{4}(h+\sqrt{s})^{-3}\big)(h')^{2}
+\tfrac{1}{2}(h+\sqrt{s})^{-3/2}|h''| \Big)e^{-(h+\sqrt{s})^{-1/2}} 
\]%
For $s\ge 1$ we use $\psi\le 2(h+\sqrt{s})$, and $\psi|\psi''| \le 3\beta((h')^{2}+|h''|)\le \beta$, while if $s\in(0,1)$
we have $\psi<2$ and 
\[%
\psi |\psi''|\le 2 |\psi''| \le 2\beta(\tfrac{3}{4}x^{5}+\tfrac{1}{4}x^{6}+\tfrac{1}{2}x^{3})e^{-x}((h')^{2}+|h''|)\le 4\beta,
\]%
where $x=(h(y)+\sqrt{s})^{-1/2}$ and we used $\beta\le 2^{-7}$, $\max_{x>0} x^{6}e^{-x}=(6/e)^{6}<125$, $\max_{x>0}
x^{5}e^{-x}=(5/e)^{5}<22$, and $\max_{x>0}x^{3}e^{-x}=(3/e)^{3}<2$.

Also, since $r<h(y)+\sqrt{\epsilon}\le h(y)+\sqrt{s}$ in $\Omega_{\epsilon}$,
\[%
\begin{aligned}%
& r^{2}|\psi''|/\psi\le 2 h^{-1}(h')^{2}+
|h''|\le 8|\big(h^{1/2}\big)'|^{2}+|h''|\le 8\beta \\ 
&r|\psi'|/\psi\le \tfrac{1}{2}h^{-1/2}|h'|=|\big(h^{1/2}\big)'|\le \beta 
\end{aligned}%
\]%
 by~\ref{h-beta}.   

Using the above inequalities in~\ref{Mv-bd} we then have 
\[%
\mathcal{M}(v_{s}) \le    \big(\!-\tfrac{\eta}{2} + (9b_{0}+4)\beta\big) \psi^{-1}\widetilde{\varphi}\hskip1pt''(r/\psi), %
\]%
and hence by taking $\beta<\eta/(36 b_{0}+16)$ we obtain
\[%
\mathcal{M}(v_{s}) \le -\tfrac{\eta}{4}\psi^{-1}\widetilde{\varphi}\hskip1pt''(r/\psi)<0,\quad s\ge \epsilon. %
\dl{super-ve}
\]%
Also by ~\ref{vtau-r} with $\widetilde{\varphi}$ in place of $\varphi$,
\[%
v_{\epsilon}-r\le \epsilon + \beta e^{-(h(y)+\sqrt{\epsilon})^{-1/2}}.  %
\dl{ved-bd}
\]%
Next we want to construct solutions $u_{t,\epsilon}$ of the SME on $\Omega_{\epsilon}$ which are $2R$-periodic in the
$y$-variable and which have boundary data given by
\[%
u_{t,\epsilon}(r,y)=\varphi_{\epsilon^{\alpha}}(r)+ %
t\bigl(v_{\epsilon}(r,y)-\varphi_{\epsilon^{\alpha}}(r)\bigr) \text{ on } \partial\Omega_{\epsilon}, %
\dl{bdry-data}
\]%
where $\alpha$ is as in~\ref{pwr-a}.

Since $\varphi_{\epsilon^{\alpha}}(r)$ is a solution of the SME when $t=0$ and $SG(\varphi_{\epsilon^{\alpha}})$ is
strictly stable, we see that the linearization of the SME at $\varphi_{\epsilon^{\alpha}}$ does not have a zero eigenvalue
with respect to zero Dirichlet data on $\partial\Omega_{\epsilon}$, so by the Implicit Function Theorem (IFT) there is a
positive solution $u_{t,\epsilon}=u_{t,\epsilon}(|\xi|,y)$ of the SME for small $t>0$ which is $2R$-periodic in $y$. In
fact notice that
\[%
\varphi_{\epsilon^{\alpha}}\le u_{t,\epsilon}\le v_{\epsilon} %
\text{ on $\Omega_{\epsilon}$ whenever $u_{t,\epsilon}$ exists,} %
\dl{ut-bds}
\]%
because $\{\varphi_{s}\}_{0<s< \epsilon^{\alpha}}$ are all solutions of the SME which are $\le u_{t,\epsilon}$ on
$\partial\Omega_{\epsilon}$, so if $u_{t,\epsilon}<\varphi_{\epsilon^{\alpha}}$ at some point of $\Omega_{\epsilon}$
then we can pick $s$ such that $\varphi_{s}\le u_{t,\epsilon}$ with equality at some point of $\Omega$, which contradicts
the maximum principle.  Similarly $\{v_{s}\}_{\epsilon \le s}$ are all supersolutions with $u_{t,\epsilon}\le v_{s}$ on
$\partial\Omega_{\epsilon}$ (by~\ref{mon-v}), so a similar maximum principle argument ensures that $u_{t,\epsilon}$
remains less than $v_{\epsilon}$ whenever $u_{t,\epsilon}$ exists.

Also, by~\ref{ved-bd} and~\ref{ut-bds},
\[%
0<u_{t,\epsilon} - r\le v_{\epsilon}-r \le \epsilon+ \beta e^{-(h(y)+\sqrt{\epsilon})^{-1/2}},\quad r<h(y)+\sqrt{\epsilon}. 
\dl{ut-r-bd}
\]%
For later reference we note that, since $u_{t_{0},\epsilon}(r,y)\ge \varphi_{\epsilon^{\alpha}}(r)>r$ and
$u_{t_{0},\epsilon}(r,y)\le v_{\epsilon}(r,y)$ (so $u_{t_{0},\epsilon}(0,r)<v_{\epsilon}(0,y)=
\epsilon+\beta e^{-(h(y)+\sqrt{\epsilon})^{-1/2}}$) by~\ref{ut-bds}, we have
\begin{align*}%
&u_{t_{0},\epsilon}(h(y)+\sqrt{\epsilon},y)/ u_{t_{0},\epsilon}(0,y)\ge (h(y)+\sqrt{\epsilon})/ u_{t_{0},\epsilon}(0,y)%
\dtg{ratio}\\ %
&\hskip0.3in  %
\ge (h(y)+\sqrt{\epsilon})/(\epsilon+\beta e^{-(h(y)+\sqrt{\epsilon})^{-1/2}}) %
\ge \tfrac{1}{2}\min\{1/\beta,1/\!\sqrt{\epsilon}\}, %
\end{align*}%
the last inequality proved by considering, for each $y$, the alternatives $\beta e^{-(h(y)+\sqrt{\epsilon})^{-1/2}}$
greater or less than $\epsilon$ and using $\min_{x>0}x^{2}e^{x^{-1}}=e^{2}/4>1$ in case of the first alternative.

Working in the set of positive $C^{2}(\widebar{\Omega}_{\epsilon})$ functions $u=u(r,y)$ ($r=|\xi|$) which are
$2R$-periodic in the $y$ variable, we use the ``method of continuity:'' For $t\in [0,1]$ and $\sigma>0$ we consider the
possibility of finding solutions $u_{t,\epsilon}$ of the SME as in \ref{bdry-data} and~\ref{ut-bds} which satisfy the
additional conditions
\[%
\left\{\begin{aligned}%
&\sup |D_{y}u_{t,\epsilon}| <\sigma \\ %
& |D(u_{t,\epsilon}(r,y)-r)| <1,\quad \tfrac{1}{2}(h(y)+\sqrt{\epsilon})\le r\le h(y)+\sqrt{\epsilon} \\ %
& \text{The linearized SME at $u_{t,\epsilon}$ has first eigenvalue positive. } %
\end{aligned}\right.%
\dl{uted-conds}
\]%
The first eigenvalue here is with respect the space of $C^{2}(\widebar{\Omega}_{\epsilon})$ functions with zero Dirichlet
data on $\partial\Omega_{\epsilon}$ and which are $2R$-periodic in the $y$-variable.

As we already mentioned, the IFT guarantees such a solution exists for all small enough $t$, so we let
$t_{0}=t_{0}(\epsilon,\beta)>0$ be defined by
\[%
t_{0}=\sup\big\{\tau\in (0,1]:u_{t,\epsilon}\text{ exists and~\ref{uted-conds} holds for all }t\in [0,\tau)\big\}. %
\dl{def-t0}
\]%
By~\cite[Example 4.1]{Sim76} (which is applicable here because of the lower and upper bounds~\ref{ut-bds}), we have a
bound on the gradient $Du_{t,\epsilon}$ for $r<\frac{1}{2}(h(y)+\sqrt{\epsilon})$ and by the second inequality
in~\ref{uted-conds} there is similarly a bound in $\frac{1}{2}(h(y)+\sqrt{\epsilon})\le r\le h(y)+\sqrt{\epsilon}$, so in fact
\[%
|Du_{t,\epsilon}(r,y)| < C, \quad (r,y)\in\widebar{\Omega}_{\epsilon}, \dl{Dute-bd}
\]%
with $C$ independent of $t$ (but possibly depending on $\epsilon$) for $t<t_{0}$.

\ref{Dute-bd} enables us to use standard local \textit{a priori} interior and boundary quasilinear elliptic estimates to bound
all the derivatives of $u_{t,\epsilon}(r,y)$ (again independent of $t$ but possibly depending on $\epsilon$).  Hence we
conclude that $u_{t_{0},\epsilon}=\lim u_{t_{k},\epsilon}$ must exist as a $C^{2}$ limit for some sequence
$t_{k}\uparrow t_{0}$.  Furthermore, if~\ref{uted-conds} holds with $t=t_{0}$ and if $t_{0}<1$, then we could use the
IFT again to show that, for some $\delta>0$, $u_{t,\epsilon}$ exists and~\ref{uted-conds} holds for all $t\in
[t_{0},t_{0}+\delta)$, contradicting the definition of $t_{0}$.  So to show $t_{0}=1$ we just have to show neither
possibility $\sup |D_{y}u_{t_{0},\epsilon}|=\sigma$ nor first eigenvalue of the linearization of the SME at
$u_{t_{0},\epsilon}$ zero nor $|Du_{t}(r,y)|=1$ can occur with $\tfrac{1}{2}(h(y)+\sqrt{\epsilon})\le r\le
h(y)+\sqrt{\epsilon}$.

We are going to prove this in case $\beta,\epsilon$ are small enough.

So let $t_{0}=t_{0}(\epsilon,\beta)\in (0,1]$ be as in~\ref{def-t0}. We claim that if $\delta>0$, if $\epsilon,\beta$ are
sufficiently small, depending on $\delta$, then for any $y_{0}\in\R$
\begin{align*}%
\dtg{c2-close}
&(1+r)^{-1}(\widetilde{u}_{t_{0},\epsilon}(r,y) - \varphi(r)) +|D(\widetilde{u}_{t_{0},\epsilon}(r,y) - \varphi(r))| \\  %
\noalign{\vskip-2pt}
&\hskip1in   + (1+r)|D^{2}(\widetilde{u}_{t_{0},\epsilon}(r,y) - \varphi(r))|   \le \delta \text{ for } |y|<1+r
\end{align*}%
and for $r<(h(y)+\sqrt{\epsilon})/u_{t_{0},\epsilon}(0,y_{0})$, where
$\widetilde{u}_{t_{0},\epsilon}(r,y)=s_{y_{0}}^{-1}u_{t_{0},\epsilon}u(s_{y_{0}} r,y_{0}+s_{y_{0}}y)$,
$s_{y_{0}}=u_{t_{0},\epsilon}(0,y_{0})$.  Geometrically this says that, in the appropriately scaled sense,
$SG(\widetilde{u}_{t_{0},\epsilon})$ is $C^{2}$ close to $SG(\varphi)$ in the region $|y|<1+r$.

Before beginning the proof of this we establish two general principles about sequences $v_{k}$ of solutions of the SME
which satisfy $v_{k}(\xi,y)>|\xi|$ on balls $B_{R_{k}}(0,0)\subset\R^{m+1}\times\R$ with $R_{k}\to \infty$ (i.e.\
$SG(v_{k}|B_{R_{k}}(0,0))$ lies in the region $|\eta|>|\xi|$ on one side of $\C_{0}\times\R=SG(r)$): We claim that for
any such sequence and for any $R_{0}>0$
\begin{align*}%
  &|(\xi_{k},y_{k})|\le R_{0} \text{ and } {\sup}_{k}v_{k}(\xi_{k},y_{k})<\infty  %
 \dtg{1st-princ}   \\  %
\noalign{\vskip-3pt}
&\hskip1in \Rightarrow {\limsup}_{k\to \infty}{\max}_{B_{R}(0,0)}v_{k}<\infty \text{ for each $R>R_{0}$,}\\ %
\noalign{\vskip-4pt \leftline{and}} %
&|(\xi_{k},y_{k})|\le R_{0} \text{ and } {\inf}_{k} (v_{k}(\xi_{k},y_{k})-|\xi_{k}|)> 0  %
\dtg{2nd-princ} \\  %
\noalign{\vskip-3pt}
&\hskip1in \Rightarrow {\liminf}_{k\to\infty}{\min}_{B_{R}} (v_{k}(\xi,y)-|\xi|) > 0 \text{ for %
  each $R>R_{0}$.} %
\end{align*}%
Both~\ref{1st-princ},~\ref{2nd-princ} are proved using the maximum principles of Solomon-White~\cite{SolW89}
and Ilmanen~\cite{Ilm96}.  For example if~\ref{1st-princ} fails then there are bounded sequences
$(\xi_{k},y_{k}),\,(\tau_{k},z_{k})$ with $v_{k}(\xi_{k},y_{k})$ bounded and $v_{k}(\tau_{k},z_{k})\to \infty$. Let
$S_{k}\subset SG(v_{k})$ be defined by $S_{k}=\{(t \xi_{k}+(1-t)\tau_{k},v_{k} (t \xi_{k}+(1-t)\tau_{k},t
y_{k}+(1-t)z_{k})\omega,ty_{k}+(1-t)z_{k}):t\in[0,1],\,\omega\in\Sph^{m}\}$.  Since $v_{k}(\tau_{k},z_{k})\to \infty$
and $R_{k}\to \infty$, we can pick $\rho_{k}\to\infty$ with $\rho_{k}^{-1}v_{k}(\tau_{k},z_{k})\to \infty$ and
$\rho_{k}^{-1} R_{k}\to\infty$, and we rescale $v_{k}$ to give
$\widetilde{v}_{k}(\xi,y)=\rho_{k}^{-1}v_{k}(\rho_{k}\xi,\rho_{k}y)$.

By~\cite[Lemma~2.3]{FouS20}, $SG(\widetilde{v}_{k})\,(=\rho_{k}^{-1}SG(v_{k}))$ has locally bounded volume in
$\R^{m+1}\times\R^{m+1}\times\R$, so the Allard compactness theorem is applicable and hence some subsequence of
$SG(\widetilde{v}_{k})$ converges in the varifold sense and locally in the Hausdorff distance sense in
$\R^{m+1}\times\R^{m+1}\times\R$ to a stationary integer multiplicity varifold $V$ with $(0,0)\in \spt
V\subset\{(\xi,\eta,y):|\eta|\ge |\xi|\}$. The maximum principle~\cite{Ilm96} ensures that $\spt V\cap
(\C_{0}\times\R)\subset\{0\}\times\R$ is impossible, so we must have $\spt V\cap
((\C_{0}\times\R)\setminus(\{0\}\times\R))\neq\emptyset$, and then, by~\cite{SolW89}, $\C_{0}\times\R\subset\spt V$. 
Let $n\in \{1,2,\ldots\}$ be the minimum $n$ such that $V$ has density $n$ on a set $T$ of positive measure in
$\C_{0}\times\R$.  But then $W=V-V_{1}$ is stationary, where $V_{1}$ is $\C_{0}\times \R$ equipped with
multiplicity identically $n$, and $V=V_{1}+W$ and hence $\density V(\xi,\eta,y)=\density V_{1}(\xi,\eta,y)+\density
W(\xi,\eta,y)$ for all $(\xi,\eta,y)$, hence $\density W(\xi,\eta,y)=0$ for $(\xi,\eta,y)\in T$. Since the density
function of a stationary varifold is upper semi-continuous, $\spt W=\{(\xi,y,\eta):\text{ density }W\ge 1\}$, and this
excludes the set of positive measure $T\subset\C_{0}\times\R$, so it is not true that $\C_{0}\times\R\subset\spt W$. 
But, $(0,0)\in \spt W$ (because the limit points of the sets $\rho_{k}^{-1}S_{k}$ include all the points
$\{(0,s\omega,0):s>0,\,\omega\in\Sph^{m}\}$, all of which are in $\spt W$), so we can repeat all of the above argument
with $W$ in place of $V$, thus concluding on the contrary that $\C_{0}\times\R\subset\spt W$.  This completes the
proof of~\ref{1st-princ}.

The proof of~\ref{2nd-princ} is similar and in fact a little simpler, because the maximum principles of~\cite{SolW89}
and~\cite{Ilm96} are applied directly to $V=\lim_{k}G(v_{k})$ without the necessity of scaling, after supposing for
contradiction that there is a bounded sequence $(\tau_{k},z_{k})$ with $v_{k}(\tau_{k},z_{k})-|\tau_{k}|\to 0$. Then the
limit points of $S_{k}=\{(t \xi_{k}+(1-t)\tau_{k},v_{k}\big(t \xi_{k}+(1-t)\tau_{k},t
y_{k}+(1-t)z_{k}\big)\omega,ty_{k}+(1-t)z_{k}):t\in[0,1],\,\omega\in\Sph^{m}\}$ for any subsequence of $k$ includes a set
$S\subset\{(\xi,\eta,y):|\eta|>|\xi|\}$ with $\,\widebar{\!S}\cap (\C_{0}\times\R)\neq \emptyset$, so we can
apply~\cite{SolW89} and~\cite{Ilm96} as in the above proof of~\ref{1st-princ} to again give a contradiction.

\smallskip

We now prove~\ref{c2-close}.  If~\ref{c2-close} fails for arbitrarily small choices of $\beta,\epsilon$ then there would be a
sequences $\beta_{k},\epsilon_{k}\downarrow 0$ and $y_{k}^{0},y_{k},r_{k}$ with
$r_{k}<(h(y_{k})+\sqrt{\epsilon_{k}})/u_{t_{0}}(0,y_{k}^{0})$, $|y_{k}|<1+r_{k}$, and
\begin{align*}%
 \dtg{c2-close-not}
&(1+r_{k})^{-1}(u_{k}(r_{k},y_{k}) - \varphi(r_{k})) +|(Du_{k})(r_{k},y_{k}) - (D\varphi)(r_{k})|  \\   %
&\hskip1.5in  + (1+r_{k})|(D^{2}{u}_{k})(r_{k},y_{k}) - (D^{2}\varphi)(r_{k}))|   \ge \delta, %
\end{align*}%
where
\[%
u_{k}(r,y)= s_{k}^{-1}u_{t_{k},\epsilon_{k}}(s_{k}r,y^{0}_{k}+s_{k}y), \quad s_{k}=u_{t_{k},\epsilon_{k}}(0,y^{0}_{k}),
\,\,t_{k}=t_{0}(\epsilon_{k},\beta_{k}). %
\dl{def-uk}
\]%
By~\ref{ratio} the domain of $u_{k}$ contains balls $B_{R_{k}}(0,0)$ with $R_{k}\to \infty$, so, since $u_{k}(0,0)=1$,
\ref{1st-princ}, \ref{2nd-princ} are both applicable and imply that $u_{k}-r$ is bounded above and below by fixed positive
constants on each ball $B_{R}(0,0)$ for all sufficiently large $k$ (depending on $R$).  Then, by~\cite[Example
4.1]{Sim76}, $|Du_{k}|$ is bounded on each ball $B_{R}(0,0)$ so a subsequence of $u_{k}$ converges locally in $C^{2}$
on all of $\R^{m+1}\times \R$ to a positive solution $u$ of the SME with $u(r,y)>r$ at every point $(r,y)$. 

Let $\rho_{k}\to\infty$. Then $\rho_{k}^{-1}u(\rho_{k} r,\rho_{k} y)>r$ and $\rho_{k}^{-1}u(0,0)\to 0$,
so~\ref{1st-princ} is again applicable, so we have local upper and lower bounds in $r>0$, and hence gradient estimates
locally in $\{(r,y):r>0\}$. Hence a subsequence of $\rho_{k}^{-1}u(\rho_{k}r,\rho_{k}y)$ converges locally smoothly in
$r>0$.  By~\ref{2nd-princ} we see that in fact the limit must be $r$ for $r>0$, otherwise there is a bounded sequence of
points $(s_{k},z_{k})$ with $\inf_{k} \rho_{k}^{-1}u(\rho_{k}s_{k},\rho_{k }z_{k})-s_{k}>0$ and then
$\rho_{k}^{-1}u(0,0)\to 0$ contradicts~\ref{2nd-princ}.

So the tangent cone of $SG(u)$ at $\infty$ is $\C_{0}\times\R$ with multiplicity~$1$ and $|D_{y}u|\le\sigma$ ($\sigma$
as in~\ref{uted-conds}), so taking $\sigma<\epsilon_{0}$, with $\epsilon_{0}$ as in the Liouville-type
theorem~\cite[Corollary 1.13]{Sim21c}, and using the fact that $u(0,0)=1=\varphi(0,0)$, we conclude that
$SG(u)=SG(\varphi)$; i.e.\ $u(r,y)=\varphi(r)$. In particular then~\ref{c2-close-not} implies $r_{k}\to\infty$. Since then
$r_{k}^{-1}u_{k}(0,0)=r_{k}^{-1}\to 0$,~\ref{2nd-princ} implies $r_{k}^{-1}u_{k}(r_{k}r,r_{k}y)-r$ converges uniformly
to $0$ in $B_{R}(0,0)$ for each $R>1$, otherwise there is a bounded sequence $(s_{k},z_{k})$ with $\inf_{k}
r_{k}^{-1}u_{k}(r_{k}s_{k},r_{k}z_{k})>0$ and then by~\ref{2nd-princ}
$\inf_{k}\inf_{B_{R}(0,0)}(r_{k}^{-1}u_{k}(r_{k}r,r_{k}y)-r)>0$, contradicting $r_{k}^{-1}u_{k}(0,0)\to 0$. In particular
we have local gradient estimates for $r_{k}^{-1}u_{k}(r_{k}r,r_{k}y)-r$ in $\{(r,y)\in B_{2}(0,0):r\ge \frac{1}{8}\}$, hence
smooth convergence to $0$ of $r_{k}^{-1}u_{k}(r_{k}r,r_{k}y)-r$ in $\{(r,y)\in B_{2}(0,0):r\ge \frac{1}{4}\}$, so
\[%
r_{k}^{-1}u_{k}(r_{k}r,r_{k}y)-r+|(Du_{k})(r_{k}r,r_{k}y)-Dr|+ r_{k}|(D^{2}u_{k})(r_{k}r,r_{k}y)-D^{2}r|\to 0
\]%
for $(r,y)\in B_{2}(0,0),\,r\ge \tfrac{1}{4}$. Since $(1,y_{k}/r_{k})\in \{(r,y)\in B_{2}(0,0):r\ge \frac{1}{2}\}$ we can take
$(r,y)=(1,y_{k}/r_{k})$ here, showing that the left side of~\ref{c2-close-not} converges to $0$, a contradiction. 
So~\ref{c2-close} is proved.

$SG(\varphi)$ is strictly stable, so~\ref{c2-close} ensures each slice $SG(u_{t_{0},\epsilon})\cap
(\R^{m+1}\times\R^{m+1}\times\{y_{0}\})$ satisfies a strict stability inequality, and then by integrating with respect to
$y_{0}$ and using the coarea formula together with the fact that $|D_{y}u_{t_{0},\epsilon}|\le \sigma$ (small) we conclude
that $SG(u_{t_{0},\epsilon})$ is strictly stable; for the details of this part of the argument we refer to~\cite[\S6]{Sim21c}.

Since we can choose $\delta$ in~\ref{c2-close} to be $<\min\{\frac{1}{2},\sigma\}$, \ref{c2-close} then also ensures that
$\sup|D_{y}u_{t_{0},\epsilon}(r,y)|\le \delta<\sigma$ for all $r<h(y)+\sqrt{\epsilon}$ and
$|D(u_{t_{0},\epsilon}(r,y)-r)|<1$ for $\frac{1}{2} (h(y)+\sqrt{\epsilon})\le r\le h(y)+\sqrt{\epsilon}$, so we have
checked~\ref{uted-conds} at $t=t_{0}$, and hence $t_{0}=1$ for all sufficiently small $\epsilon,\beta$ as claimed.

Because $|D\varphi|<1$,~\ref{c2-close} (with $t_{0}=1$) also implies
\[%
{\sup}_{0\le r\le h(y)+\sqrt{\epsilon}}|Du_{1,\epsilon}|<1+\delta<1+\sigma. %
\dl{grad-est}
\]%
Because of~\ref{grad-est} and~\ref{ut-r-bd} we can apply standard local quasilinear elliptic interior and boundary estimates
to bound the derivatives of the difference $u_{1,\epsilon}(r,y)-r$ in regions $\theta (h(y)+\sqrt{\epsilon})\le r\le
h(y)+\sqrt{\epsilon}$; preparatory to this observe that if $y\in \R$ and $t\in [-4,4]$ then
$h(y+t(h(y+\sqrt{\epsilon})))-h(y)=t h'(y+s(h(y+\sqrt{\epsilon})))(h(y)+\sqrt{\epsilon})$ for some
$s$ between $0$ and $t$, so, for $|y-z|\le 4 (h(y)+\sqrt{\epsilon})$,
\[%
\tfrac{1}{2} (h(y)+\sqrt{\epsilon})\le(1-4\beta) (h(y)+\sqrt{\epsilon})\le h(z)+\sqrt{\epsilon}\le (1+4\beta) %
(h(y)+\sqrt{\epsilon})\le 2(h(y)+\sqrt{\epsilon}). %
\dl{h-comp} 
\]%
So in particular, using the right hand inequality, 
\[%
e^{-(h(z)+\sqrt{\epsilon})^{-1/2}} \le  e^{-\frac{1}{2}(h(y)+\sqrt{\epsilon})^{-1/2}}, \quad  |y-z|\le 4
(h(y)+\sqrt{\epsilon}),
\dl{bd-E}
\]%
and using \ref{grad-est}, \ref{ut-r-bd}, \ref{bd-E} and standard elliptic estimates for the difference $u_{1,\epsilon}-r$, we
have, for any $\theta\in \big(0,\tfrac{1}{2}\big]$,
\[%
\big|D^{k}\big(u_{1,\epsilon}(r,y)-r\big)\big|\le C_{k}\theta^{-k}(h(y)+\sqrt{\epsilon})^{-k} %
\big(\epsilon+\beta e^{-\frac{1}{2}(h(y)+\sqrt{\epsilon})^{-1/2}}\big) %
\dl{deriv-est}
\]%
for $\theta (h(y)+\sqrt{\epsilon}) \le r \le h(y)+\sqrt{\epsilon}$ and $k=1,2,3,\ldots$, with $C_{k}$ independent of
$\theta$.

Up to this point all the constructions have been made assuming the function $h=h_{R}$ and $K=K_{R}$ are
$2R$-periodic, where $R>1$ was fixed, and correspondingly then $u_{1,\epsilon}=u_{1,\epsilon,R}$ also depends on $R$.
 Suppose now that we have an arbitrary closed $K\subset\R$ and arbitrary non-negative $C^{\infty}$ $h$
satisfying~\ref{h-beta} with $\beta/3$ in place of $\beta$ and with $K=\{y:h(y)=0\}$.  Let $\zeta:\R\to[0,1]$ be
$C^{\infty}$ with $\zeta(y)=0$ for $|y|\ge \frac{7}{8}$, $\zeta(y)>0$ for $|y|<\frac{7}{8}$, $\zeta(y)=1$ for
$|y|<\frac{1}{2}$, and then let $\zeta_{R}(y)=\zeta(y/R)$, and $h_{R}=$ the $2R$-periodic extension of
$\zeta_{R}^{2}h\big|[-R,R]$.  Then, for all $R$ sufficiently large,~\ref{h-beta} holds for $h_{R}$ (possibly with new
constants $C_{k}$, but $C_{k}$ independent of $R$), and so the constants $C_{k}$ in~\ref{deriv-est} can also be chosen
independent of $R$.

So let $\epsilon=\epsilon_{k}\to 0$ and $R=R_{k}\to \infty$. Since the gradient estimate \ref{grad-est} and the derivative
estimates in~\ref{deriv-est} are independent of $\epsilon$ and $R$, some subsequence of $u_{1,\epsilon_{k},R_{k}}$
converges uniformly in $\widebar{\Omega}$ and locally smoothly in $\widebar{\Omega}\setminus(\{0\}\times\R)$ to
some $u_{1}$ which has bounded gradient in $\Omega$ and is smooth in $\Omega\setminus(\{0\}\times\R)$ with
$u_{1}(r,y)\ge r$; but by the Hopf maximum principle applied to the difference $u_{1}-r$ (which is the difference of two
solutions of the SME) we see that strict inequality $u_{1}(r,y)>r$ must hold in $\Omega\setminus(\{0\}\times\R)$,
because by construction $u_{1}(r,y)=\smash{\beta e^{-h^{-1/2}(y)} \widetilde{\varphi}(h(y)/\beta e^{-h^{-1/2}})}>0$ on
$\partial\Omega$.  $SG(u_{1})$ is stationary and contained in $\{(\xi,\eta,y):|\eta|>|\xi|\}\cup( \{0\}\times(\R\setminus
K))$, so the maximum principle~\cite{Ilm96} implies $SG(u_{1})\cap (\C_{0}\times(\R\setminus K))=\emptyset$, so in
fact $u_{1,\epsilon_{k}}$ is bounded below by a fixed positive constant on each compact subset $\Gamma\subset
\Omega$, hence $u_{1,\epsilon_{k},R_{k}}$ to $u_{1}$ locally in $C^{\infty}(\Omega)$, and $u_{1}$ is positive and
smooth on all of $\Omega$ with
\[%
\left\{\begin{aligned}%
&\mathcal{M}(u_{1})(r,y)=0, \,\,|Du_{1}(r,y)|<1+\sigma,\,\,\, (r,y)\in \Omega \\ %
&0< u_{1}(r,y)-r \le \beta e^{-h^{-1/2}(y)},\,\,\,(r,y)\in \Omega %
\,\,\text{ (by~\ref{ut-r-bd})} \\ %
&u_{1}(h(y),y)=\beta e^{-h(y)^{-1/2}}\widetilde{\varphi}\big(h(y)/\beta e^{-h(y)^{-1/2}}\big),
\,\,y\in\R\setminus K.
\end{aligned}\right.%
\dl{props-u10}
\]%
Also, by~\ref{deriv-est} with $\theta=\frac{1}{4}$, 
\[%
\big|D^{k}\big(u_{1}(r,y)-r\big)\big| \le C_{k,\ell}\beta h(y)^{\ell},
 \quad \tfrac{1}{4}h(y)\le r\le h(y), \,\, k,\ell=1,2,\ldots. %
\dl{props-Du10}
\]%
Now take a $C^{\infty}$ function $\zeta:\R\to [0,1]$ with $\zeta(t)=0$ for $t\le \frac{1}{2}$,
$0<\zeta(t)<1,\,t\in (\frac{1}{2},1)$ and $\zeta(t)=1$ for $t\ge 1$, and for $(r,y)\in [0,\infty)\times (\R\setminus K)$ 
define
\[%
u(r,y) =\zeta(\vphantom{a^{2}}r/h(y))r + (1-\zeta(\vphantom{a^{2}}r/h(y)))u_{1}(r,y) %
= u_{1}(r,y) -\zeta(\vphantom{a^{2}}r/h(y))\bigl(u_{1}(r,y)-r\bigr).
\]%
Thus $u(r,y)-r = \left(\vphantom{a^{2}}u_{1}(r,y)-r\right)\,\left(1-\zeta(\vphantom{a^{2}}r/h(y))\right)$ ($=u_{1}$ for
$r<\frac{1}{2}h(y)$) and hence by~\ref{props-u10} and~\ref{props-Du10} $u(r,y)$, for $y\in \R\setminus K$, has the
properties (with new constants $C_{k,\ell}$)
\[%
\left\{\begin{aligned}%
&\!\mathcal{M}(u)(r,y)=0,\,\, r<\tfrac{1}{2}h(y),\,\,\,\,  0< u(r,y)-r \le \beta e^{-h^{-1/2}(y)},\,\,  r<h(y),  \\ %
&\,u(r,y)=r, \,\,\,r\ge h(y)  \\ %
&\!\left|D_{r,y}^{k}\left(u-r\right)\right| \le C_{k,\ell}\beta h(y)^{\ell}, \,\,\,\, %
\tfrac{1}{4}h(y) < r ,\,\, k,\ell=1,2,\ldots  %
\end{aligned}\right.%
\dl{props-u}
\]%
Also the above estimates enable us to extend $u$ to be $C^{\infty}$ on all of
$(\R^{m+1}\times\R^{m+1}\times\R)\setminus(\{0\}\times\{0\}\times K)$ by taking
\[%
u(r,y) = r, \quad y\in K, 
\dl{ext-u}
\]%
and then clearly $\sing G(u)=\{0\}\times\{0\}\times K$. 

Finally we remark that the slicing argument used to check strict stability of $SG(u_{1,\epsilon})$ (hence $SG(u_{1})$)
extends without change to $SG(u)$, because for each $y_{0}\in \R\setminus K$, $u$ evidently satisfies an inequality like
the inequality for $u_{1,\epsilon}$ in~\ref{c2-close} for all $r\ge 0$.

\section{Construction of a Metric for $\R^{m+1}\times\R^{m+1}\times\R$}

With $u$ as in~\ref{props-u} and \ref{ext-u}, we can now discuss how to find a metric $g$ with respect to which $SG(u)$
is minimal.

We look for a metric $g$ of the form
\[%
g_{|(\xi,\eta,y)}={\sum}_{i=1}^{m+1}d\xi_{i}^{2}+f(\xi,y){\sum}_{j=1}^{m+1}d\eta_{j}^{2}+dy^{2}, %
\dl{g-f}
\]%
where $f$ is smooth, close to $1$ everywhere, and equal to $1$ on $\R\times K$; also we shall construct $f$ so that it is a
function of $|\xi|,y$, so henceforth we consider only functions $f=f(r,y)$, $r=|\xi|$.  Applying the area formula, except that
now we use the metric $g$ of~\ref{g-f} rather than the standard Euclidean metric, the area functional is then
\[%
\mu_{g}(SG(v))=\sigma_{m}\int_{V}\sqrt{1+f\,|Dv|^{2}}\,\,f^{m/2}v^{m}\,d\xi dy %
\dl{f-area}
\]%
for any positive $C^{2}$ function $v$ on a domain $V\subset\R^{m+1}\times\R$, where $\mu_{g}$ denotes
$(2m+2)$-dimensional Hausdorff measure with respect to the metric $g$ for $\R^{m+1}\times\R^{m+1}\times\R$ and
$\sigma_{m}$ is the volume of the unit sphere $\Sph^{m}$ in $\R^{m+1}$.  Thus the Euler-Lagrange equation for the
functional on the right of~\ref{f-area} is equivalent to the statement that the symmetric graph $SG(v)$ is a minimal (zero
mean curvature) hypersurface relative to the metric $g$ for $\R^{m+1}\times\R^{m+1}\times\R$.

By direct computation, the Euler-Lagrange equation is
\[%
\tfrac{1}{2} \big(m+1+\frac{1}{1+f\,|Dv|^{2}}\big)\,Df\cdot Dv = %
-f\big(\Delta v-f\frac{v_{r}^{2}v_{rr}+v_{y}^{2}v_{yy}+2v_{r}v_{y}v_{ry}}{1+f|Dv|^{2}}\big)+\frac{m}{v}. %
\]%
At points where $Dv\neq 0$ 
this is a non-degenerate first-order PDE for $f$, and we can impose the initial condition
\[%
f=1 \text{ on the hypersurface } |\xi|=\tfrac{1}{2}h(y), %
\dl{init-f}
\]%
i.e.\ $f=1$ on inner boundary of the transition region.

We want to show that we can choose $f$ so that this holds with the choice $v=u$, $u$ as in~\ref{props-u}, \ref{ext-u}. 
With such a choice of $v$ and with $\,\,\widebar{\!\!f}=1-f$, the equation for $f$ can be rewritten
\begin{align*}%
&\,\,\,\tfrac{1}{2}\!\big(m+1+\frac{1}{1+|Du|^{2}-\,\,\vphantom{\bigl(}\swidebar{\!\!f}\,|Du|^{2}}\big) %
Du\cdot D\,\,\widebar{\!\!f}  ={\cal{}M}(u) -E(u,\,\,\widebar{\!\!f})\,\,\widebar{\!\!f},\\
\noalign{\vskip3pt} &E(u,\,\,\widebar{\!\!f})\!=\!\Delta_{\R^{m+1}\times\R}u-\tfrac{2+|Du|^{2}- %
  \,\,\swidebar{\!\!f}(1+|Du|^{2})} {(1+|Du|^{2}-\,\,\swidebar{\!\!f}|Du|^{2}) %
  (\vphantom{\widebar{\!\!f}}1+|Du|^{2})} %
(u_{r}^{2}u_{rr}+u_{y}^{2}u_{yy}+2u_{r}u_{y}u_{ry}).
\end{align*}%
The inhomogeneous term in this equation is $\mathcal{M}(u)$ which is zero in $|\xi|<\frac{1}{2} h(y)$ by~\ref{props-u},
so we must have $\,\,\widebar{\!\!f}$ identically zero on $\frac{1}{4}h(y)\le |\xi|\le \frac{1}{2} h(y)$ (by the uniqueness
theorem for quasilinear first order PDE's); also, since $E(r,\,\,\widebar{\!\!f})=m /r$ and $\mathcal{M}(r)=0$ ($r=|\xi|$),
the above implies
\[%
\tfrac{1}{2}\big(m+1+\frac{1}{1+|Du|^{2}(1-\,\,\vphantom{\bigl(}\overline{\!\!f}\,)}\big)\,Du\cdot D\,\,\widebar{\!\!f}%
 =-\tfrac{m}{r}\,\,\widebar{\!\!f} +\left(\vphantom{\widebar{\!\!f}}{\cal{}M}(u)-\mathcal{M}(r)\right) %
 -\left(E(u,\,\,\widebar{\!\!f})-E(r,\,\,\widebar{\!\!f})\right)\,\,\widebar{\!\!f}, %
\]%
and, since $Du=Dr+D(u-r)=r^{-1}\xi+a$ with $a=D(u-r)$, by~\ref{props-u} this can be written
\[%
\tfrac{1}{2}\big(m+1+\frac{1}{1+\left(1+|a|^{2}-2a\cdot \xi/r\right)(1-\,\,\vphantom{\bigl(}\overline{\!\!f}\,)}\big)\,  %
 \left(r^{-1}\xi+a(\xi,y)\right)\cdot D\,\,\widebar{\!\!f}   %
 =-m\,\,\widebar{\!\!f}/r+b(\xi,y,\,\,\widebar{\!\!f})\,\,\,\widebar{\!\!f} + c(\xi,y,\,\,\widebar{\!\!f})
 \dl{f-pde}  %
\]%
with
\[%
\bigl|a(\xi,y)\bigr|\le \beta\, (<\tfrac{1}{2}), \quad r\ge \tfrac{1}{4}h(y), \dl{non-deg}
\]%
provided we take $\beta$ in~\ref{h-beta} and~\ref{props-u} small enough, thus ensuring that~\ref{f-pde} is uniformly
non-degenerate in the region $r>\frac{1}{4}h(y)$, and, again by~\ref{props-u},
\[%
\bigl|D_{\xi,y}^{k}a(\xi,y)\bigr|+{\sup}_{|z|\le \frac{1}{2}}\big(\bigl|D_{\xi,y,z}^{k}b(\xi,y,z)\bigr| +
\bigl|D_{\xi,y,z}^{k}c(\xi,y,z)\bigr|\big) \le C_{k,\ell}\beta h(y)^{\ell}, \dl{bds-ab}
\]%
for $r\ge \tfrac{1}{4}h(y)$ and $k,\ell=0,1,2,\ldots$, with constants $C_{k,\ell}$ independent of $\beta$.

Then using~\ref{non-deg} and~\ref{bds-ab} together with standard estimates for non-degenerate first order PDE's (as
discussed in~\cite{Sim23}), we have
\[%
\bigl|D_{r,y}^{k}\,\,\widebar{\!\!f}(r,y)\bigr|\le C_{k,\ell}\beta h(y)^{\ell},  \,\,\,k,\ell=0,1,2,\ldots, %
\dl{bds-df}
\]%
$\tfrac{1}{4}h(y)\le r\le h(y)$, with suitable constants $C_{k,\ell}$ depending only on $k,\ell$ (and not depending on
$\beta$).

Also, since $u=r$ for $r\ge h(y)$ (so $a=b=c=0$ for $r\ge h(y)$), the PDE for $\,\,\widebar{\!\!f}$ can be written for $r\ge
h(y)$ as the ODE
\[%
\left(m+1+(2-\,\,\widebar{\!\!f})^{-1}\right)\,\,\widebar{\!\!f}_{r} = -2m \,\,\widebar{\!\!f}/r,  %
\]%
which can be integrated to give
\[%
\left(2- \,\,\widebar{\!\!f}(r,y)\right)^{-\kappa_{1}}\, %
\widebar{\!\!f}(r,y) =\left(2- \,\,\widebar{\!\!f}(h(y),y)\right)^{-\kappa_{1}} %
\widebar{\!\!f}(h(y),y)\left(\vphantom{\widebar{\!\!f}}h(y)/r\right)^{\kappa_{2}}\!\!, \,\, r\ge h(y), %
\]%
where $\kappa_{1} =\frac{1}{2m+3}$, $\kappa_{2}= \frac{4m}{2m+3}$, and hence the inequalities~\ref{bds-df} are also
valid (with possibly different constants $C_{k,\ell}$) for $r>h(y)$.

Thus \ref{bds-df} holds for $r>\frac{1}{4}h(y)$, and since $f$ is identically $1$ on $\tfrac{1}{4}h(y)\le r\le
\tfrac{1}{2}h(y)$, we can extend $f$ to be identically $1$ in the region $0\le r<\tfrac{1}{2}h(y)$.  So $f$ exists and is
smooth and positive on $\R^{m+1}\times (\R\setminus K)$ and $|D_{r,y}^{k}\,\,\widebar{\!\!f}(r,y)|\le C_{k,\ell}\beta
h(y)^{\ell}$, $(r,y)\in [0,\infty)\times(\R\setminus K)$, $k,\ell=0,1,2,\ldots$, and so we can extend $f$ to be $C^{\infty}$
on all of $\R^{m+1}$ by taking $f=1$ on $\R^{m+1}\times K$.

\providecommand{\bysame}{\leavevmode\hbox to3em{\hrulefill}\thinspace}

\vskip.2in

{\obeylines \multiply \baselineskip by 11\divide \baselineskip by 14 
Mathematics Department, Stanford University 
Stanford CA 94305, USA 
{\small lsimon@stanford.edu}}

\end{document}